\documentclass[12pt]{amsart}
\usepackage{amssymb,array,epsfig}

\newtheorem{lemma}{Lemma}
\newtheorem{proposition}[lemma]{Proposition}
\newtheorem{theorem}[lemma]{Theorem}

\newcommand{\CE}{\mathcal{E}}
\newcommand{\CH}{\mathcal{H}}
\newcommand{\CL}{\mathcal{L}}

\newcommand{\R}{\mathbb{R}}
\newcommand{\C}{\mathbb{C}}

\newcommand{\N}{\mathbb{N}}

\newcommand{\spec}{\mathrm{Spec}\,}
\newcommand{\spe}{\mathrm{Spec}}

\newcommand{\num}{\mathrm{Num}\,}
\newcommand{\dom}{\mathrm{Dom}\,}
\newcommand{\ran}{\mathrm{Ran}\,}

\renewcommand{\Re}{\mathrm{Re}\,}
\renewcommand{\Im}{\mathrm{Im}\,}
\newcommand{\ess}{\mathrm{ess}}
\newcommand{\disc}{\mathrm{disc}}
\newcommand{\dist}{\mathrm{dist}\,}
\newcommand{\ud}{\,\mathrm{d}}

\title[Non-variational...]{Non-variational
approximation of discrete eigenvalues of self-adjoint operators}

\author[L. Boulton]{Lyonell Boulton$^{1}$}

\date{March 2005}

\thanks{$^1$Pacific Institute for the Mathematical Sciences
postdoctoral fellow.}



\begin{document}

\begin{abstract}
We establish sufficiency conditions in order to achieve
approximation to discrete eigenvalues of self-adjoint operators in
the second-order projection method suggested recently by Levitin and
Shargorodsky, \cite{lesh}.  We find explicit estimates for the
eigenvalue error and study in detail two concrete model examples.
Our results show that, unlike the majority of the standard methods,
second-order projection strategies combine non-pollution and
approximation at a very high level of generality.
\end{abstract}


\maketitle

\

\vspace{-.5in}

\section{Introduction}

Let $M$ be a self-adjoint operator acting on a dense domain, $\dom
M$, in a separable Hilbert space $\CH$. Let the spectrum of $M$ be
\linebreak split into isolated eigenvalues of finite multiplicity
(discrete spectrum) and degenerate points (essential spectrum), in
symbols \linebreak $\spec M= \spe_\disc M \cup \spe_\ess M.$ In
order to approximate discrete eigenvalues $\lambda$ of $M$, we may
consider using the following projection method: choose an
orthonormal basis of $\CH$, $\{\phi_k\}_{k=1}^\infty \subset \dom
M$, and find the spectrum of large matrices $M_n$ resulting from
compressing $M$ to the finite dimensional subspaces
$\CL_n=\mathrm{span}\{\phi_1,\ldots,\phi_n\}$.

The successful outcomes of this strategy is illustrated by the well
known Rayleigh-Ritz theorem for the approximation of variational
eigenvalues, see e.g. \cite[Theorem~XIII.4]{resi}. Assume that $M$
has a non-degenerate ground eigenvalue $\lambda =\min [\spec M] <
\min [\spec\!_{\mathrm{ess}}\, M]$. Let $\Pi_n$ be the orthogonal
projection onto $\CL_n$ and $M_n:=\Pi_nM|\CL_n$. If
\begin{equation} \label{e5}
\lim_{n\to \infty}\|M_n\Pi_n \psi-\lambda \psi\| = 0 \qquad
\mbox{for\ all}\ M\psi=\lambda \psi,
\end{equation}
then the first eigenvalue of $M_n$ converges from above to $\lambda$
(i.e. approximation is achieved) and the second eigenvalue of $M_n$,
counting multiplicity, can not be smaller than  $\min[\spec M
\setminus \{\lambda\}]$ (i.e. no chance of spurious eigenvalues).
Conditions such as \eqref{e5} are useful in applications because
they can sometimes be verified on abstract grounds,
\cite[Ex.~XIII.2.1 and XIII.2.2]{resi}. Similar results hold for any
other variational eigenvalue (i.e. those below the minimum of the
essential spectrum or above the maximum of this set) simple or
otherwise.

Although this technique is proven to be valuable in the study of
spectral properties of self-adjoint operators, think e.g. of
Galerkin approximation, two main difficulties are in place when we
extend it to the analysis of non-variational eigenvalues:
\begin{itemize}
\item[-] \emph{Lack of approximation:} no
eigenvalue of $M_n$ is close to $\lambda$ in the limit $n\to\infty$,
and
\item[-] \emph{Spectral pollution:}  there are eigenvalues
$z_n$ of $M_n$ which appear to converge to some $\mu$ not in the
spectrum of $M$.
\end{itemize}
In fact any $\mu \in \mathrm{conv} [\spe_\ess M] \setminus \spec M$
is potentially an accumulation point of $\spec M_n$ in the limit
$n\to \infty$. This is known to occur in several important
applications such as elasticity, electromagnetism and hydrodynamics,
see e.g. \cite{arn}, \cite{boffi1}, \cite{sp1}, \cite{sp2} and the
references therein. Widely available commercial packages such as
FEMLAB are known to produce spectacularly incorrect results even in
cases of the simplest one-dimensional Stokes' type systems,
\cite{lesh}. Actually one may construct simple examples of bounded
$M$ and natural orthonormal basis $\{\phi_k\}$ with the property
that \emph{every} $\mu\in \mathrm{conv} [\spe_\ess M] \setminus
\spec M$ is an accumulation point of $\spec M_n$ in the limit $n\to
\infty$, see \cite{lesh} and \cite{bou2}.

\medskip

The standard approach to deal with pollution around  non-variational
eigenvalues, aims at choosing $\CL_n$ sufficiently close to the
spectral subspaces of $M$ in order to capture the eigenfunctions
associated to $\lambda$, see e.g. \cite{arn}, \cite{boffi1} and
\cite{sp2}. Unfortunately no universal device is known for
constructing such subspaces.

A different approach hinges on the spectral theorem and it exploits
the fact that for $\zeta$ lying in a gap of
$\spec\!_\mathrm{ess}\,M$, the eigenvalues of $\zeta-M$ inside the
corresponding shifted gap enclosing the origin, become variational
eigenvalues of $(\zeta-M)^2$. In \cite{lesh} Levitin and
Shargorodsky consider this idea in a concrete manner. They propose
complementing the usual projection method with \emph{a posteriori}
estimates found by computing the eigenvalues of the matrix
polynomial
\begin{equation} \label{e10}
  P_n(z):=\Pi_n(z-M)^2|{\CL_n}=z^2-2M_nz+[M^2]_n.
\end{equation}
Usually $\spec P_n$ contains non-real points (see the definitions in
Section~2). Their result is based on a remarkable property
\cite[Theorem~2.6]{lesh} guaranteing that for $z\in \spec P_n$,
\begin{equation} \label{e31}
[\Re z\!\!-\!\! |\Im z|,\Re z\!\!+\!\!|\Im z|] \cap \spec M \not =
\varnothing .
\end{equation}
Hence the strategy would be finding points in $\spec P_n$. If these
points are close to $\R$, they must also be close to the spectrum of
$M$.

Numerical evidence found in \cite{bou1} and \cite{lesh} suggests
that approximation occurs for concrete operators and natural basis
$\{\phi_k\}$. The aim of the present paper is to discuss in detail
this complementary problem of approximation in this ``second-order''
projection method. Our main contribution summarizes in the following
statement: \emph{if $\lambda\in \spec\!_\mathrm{disc}\, M$ and
\begin{equation} \label{e6}
  \lim_{n\to \infty}\|\Pi_n M^k\Pi_n\psi-\lambda^k\psi\|
= 0,
  \quad k=0,1,2,\ {\it whenever}\ M\psi=\lambda \psi,
\end{equation}
then there exists $z_n\in\spec P_n$ such that $z_n\to \lambda$ as
$n\to \infty$}. See Theorems~\ref{t1} and \ref{t5} below and also
\cite{bou2}. Hence, discrete eigenvalues satisfying \eqref{e6} will
always be approached no matter their location relative to
$\spec\!_\mathrm{ess} M$. Notice that condition \eqref{e6} here
plays the role of \eqref{e5} in the ``linear'' projection method.

\medskip

In Section~2 we establish the main theoretical contribution.
Theorem~\ref{t1} finds explicit estimates for the eigenvalue error
$|z_n-\lambda|$ in terms of bounds for the difference in the left
side of \eqref{e6}. This result depends on knowing exactly the
entries of the matrices $M_n$ and $[M^2]_n$, hence its scope is
mainly theoretical. Theorem~\ref{t5} on the other hand, demonstrates
that introducing noise to the entries of $P_n(z)$ within a certain
tolerance, does not affect approximation. This result is better
suited for applications.  We give precise bounds for this tolerance
in terms of invariants of the problem.

Section~3 illustrates the scope of the theoretical results by means
of two model examples. The first one corresponds to finite rank
perturbations of multiplication operators by a bounded symbol. We
find explicit estimates for the error $|z_n- \lambda|$ when
$\{\phi_k\}$ is of Fourier type and report on numerical outputs. The
second model is unbounded, the Schr\"odinger operator with band gap
essential spectrum. We choose $\{\phi_k\}$ to be the Hermite
functions, then find explicit upper bounds for eigenvalue error in
terms of smoothness properties of the potential.

In order to make the paper more readable, the proofs of
Theorems~\ref{t1} and \ref{t5} are to be found in Section~5 while
Section~4 introduces and describes a convenient notation in order to
simplify the presentation of most arguments. Being of obvious
interest, in Section~6 we discuss the question of whether
approximation also occurs for the essential spectrum.


\section{Speed of convergence and stability: general results}

We begin this section by fixing some notation. Let $\lambda \in
\spec_\mathrm{disc} M$. Throughout this paper $\{\Pi_n\}$ denotes a
family of orthogonal projections onto $\CH$ satisfying \eqref{e6}
such that $\CL_n:=\ran \Pi_n \subset \dom M^2$. The sequence
$\delta(n)>0$ shall always denote an upper bound for the left side
of \eqref{e6} such that $\delta(n)\to 0$ as $n\to \infty$: for all
$\psi \in \dom M$ satisfying $M\psi=\lambda \psi$ and $\|\psi\|=1$,
\begin{equation} \label{e2}
  \|\Pi_nM^k\Pi_n\psi-\lambda^k\psi\|\leq \delta(n)\qquad
  \mathrm{for}\ n\in \N\ \mathrm{and}\
   k=0,1,2.
\end{equation}
Below we shall always assume that $\psi$ is an eigenvector
associated to $\lambda$ normalized by $\|\psi\|=1$.

Let $P(z)=\sum_{k=0}^m A_k z^k$ where $z$ is a complex variable,
$A_k\in \C^{n\times n}$ and $\det A_m\not =0$. The spectrum of $P$
is, by definition, the set of eigenvalues
\[
 \spec P:=\{z\in \C:\det P(z)=0\}.
\]
The hypothesis $\det A_m\not=0$ ensures that $\spec P$ comprises no
more than $mn$ finite points.

We are concerned with the spectrum of the quadratic matrix
polynomial $P_n(z)$ given by \eqref{e10}. The determinant of
$P_n(z)$ is a scalar polynomial in $z$ of degree $(2\dim \CL_n)$, so
that $\spec P_n$ is a finite set comprising at most $(2\dim \CL_n)$
different complex numbers. In general these points do not intersect
the real line, unless $\CL_n$ contains an eigenfunction of $M$.
Furthermore, since $\overline{\det P_n(z)} = \det
P_n(\overline{z}),$ $\R$ is an axis of symmetry for $\spec P_n$. In
Section~4 we characterize $\spec P_n$ as the poles of a certain
positive-valued subharmonic function. Other descriptions of the
spectrum of matrix polynomials  better suited for computation, e.g.
as the eigenvalues of block matrices, may be found in the
literature, see \cite{glr} or \cite[Part II]{glr2}.

The goal of the procedure suggested by Levitin and Shargorodsky in
\cite{lesh} is finding the points in $\spec P_n$ which are close to
$\R$. Statement \eqref{e31} ensures that these points will
necessarily be close to $\spec M$. A complementary assertion
guarantees that discrete eigenvalues satisfying \eqref{e6} will
always be approached no matter their location relative to
$\spec\!_\mathrm{ess} M$. Our first theorem in this direction
establishes that the rate of convergence $|z_n-\lambda|\to 0$ for
$z_n\in \spec P_n$ is at least a power $1/2$ the rate at which the
$\CL_n$ approximate $M^k\psi$, $k=0,1,2$. It is unclear whether the
claimed power $1/2$ is sharp for general choices of $\CL_n$. In
Section~\ref{ss2} below, we provide specific numerical evidence
suggesting that for particular cases the actual eigenvalue error is
$|z_n-\lambda|=O(\delta(n)^\alpha)$ for $\alpha\approx 1$.

\begin{theorem} \label{t1}
Let $\lambda \in \spec_\disc M$ and assume that \eqref{e2} holds.
There exists $b>0$ independent of $n$ and $z_n\in\spec P_n$, such
that
\begin{equation} \label{e20}
   |z_n-\lambda|< b [\delta(n)]^{1/2}, \qquad n\in \N.
\end{equation}
\end{theorem}

\medskip

The scope of Theorem~\ref{t1} is mainly theoretical. Computing
$\spec P_n$ usually requires estimating the coefficients of
$P_n(z)$. Since $P_n(z)$ is Hermitian for all $z\in\R$, a well known
result  in the theory of matrix polynomials (cf.
\cite[Theorem~II.2.6]{glr2}), ensures that we can always find a
factorization
\begin{equation}\label{e12}
P_n(z)=(z-A_n)(z-A_n^\ast),\qquad z\in \C, \end{equation} where
$A_n$ are square matrices of size $\dim \CL_n$ constant in $z$.
Clearly \linebreak $\spec P_n=(\spec A_n) \cup \overline{(\spec
A_n)}$. A small perturbation in the coefficients of $P_n$ will
change the coefficients of $A_n$. Typically the eigenvalues of $A_n$
are not semi-simple. In fact their condition number might, in
principle, be large as $n$ is large, forcing $\spec P_n$ to be
sensitive to small changes in the entries of $M_n$ and $[M^2]_n$,
see \cite{gol}. Our next result establishes that this sensitivity
can be controlled uniformly in neighbourhoods of
$\spe_{\mathrm{disc}} M$. To be more specific, approximation to a
small $\delta$-neigh\-bourhood of $\lambda$ is always achieved, if
we estimate the coefficients of $P_n$ within an error smaller than
some prescribed tolerance, $w_k\varepsilon$.

Below and elsewhere we shall write $\mu:=\dist [\lambda,\spec
M\setminus \{\lambda\}]>0$. The norm $\|\!\cdot\!\|$ for matrices
shall always refer to the uniform operator norm. Let
$\varepsilon\geq 0$ and $\bar{w}=(w_0,w_1)$ where $w_k \geq 0$. We
denote by $\mathcal{P}_{\varepsilon,\bar{w}}$ the set of sequences
of linear matrix polynomials $(Q_n)_{n=1}^\infty$ such that
$Q_n(z)=F_nz-G_n$, where $F_n$ and $G_n$ are square matrices
constant in $z$ of size $(\dim \CL_n)$ satisfying $\|G_n\|\leq
w_0\varepsilon$ and $\|F_n\|\leq w_1\varepsilon$.

\begin{theorem} \label{t5}
Let $\lambda \in \spec_\disc M$. Assume that \eqref{e6} holds and
that $w_0,w_1\geq 0$ do not vanish simultaneously. Let
$0<\delta<\mu/4$ be fixed.  Let
\begin{equation} \label{e7}
   0\leq \varepsilon <\frac{\delta^2 \mu^2}{2(2\delta^2+3\mu^2)
   [w_0+w_1(\mu/4+|\lambda|)]}.
\end{equation}
There exists $N>0$ only dependant upon $\delta$, $w_0$ and $w_1$,
such that for all $(Q_n)_{n=1}^\infty \in
\mathcal{P}_{\varepsilon,\bar{w}}$,
\begin{gather*}
\spec (P_n+Q_n)\cap
\{\delta \leq |z-\lambda|\leq
\mu/4\}=\varnothing \qquad \mathrm{and}\\
\spec (P_n+Q_n)\cap \{|z-\lambda|<\delta\}\not = \varnothing \qquad
\mathrm{for}\ n>N.
\end{gather*}
Moreover, if we count multiplicities, the number of eigenvalues of
$P_n$ and $P_n+Q_n$ in $\{|z-\lambda|<\delta\}$ coincide.
\end{theorem}

In other words, if we aim at detecting $\lambda$ with an error of
$\delta$, it is enough to consider estimations
\[
 P_n(z)+Q_n(z)=z^2-(2M_n-F_n)z+([M^2]_n-G_n)
\]
with sufficiently small $\|F_n\|$ and $\|G_n\|$, and find the
spectrum of $P_n+Q_n$ for sufficiently large $n$. The weights $w_0$
and $w_1$ are introduced in order to allow independent control on
how perturbations are measured in the truncations of $M$ and $M^2$.
For instance, two possibilities are: the absolute weights
$w_0=w_1=1$ and the relative weights
\[
 w_0=\frac{\mu^2}{4(2\delta^2+3\mu^2)}
 \qquad \mathrm{and} \qquad w_1=\frac{w_0}{(\mu/4+|\lambda|)}.
\]

\section{Speed of convergence: examples}

In the examples presented below we find $\delta(n)$ explicitly for
concrete operators $M$. Among other results, we illustrate how the
error $|z_n-\lambda|$ depends strongly on the correct choice of
$\CL_n$. We also provide numerical evidence suggesting that, in
cases, the power $1/2$ of $\delta(n)$ predicted by Theorem~\ref{t1}
can actually be improved to a power $\alpha\approx 1$.

\subsection{Finite rank perturbations of multiplication operators}
\label{ss1} As a first example we consider $M=S+K$ acting on $f\in
L^2(-\pi,\pi)=:L^2$, where
\[
   Sf(x)=s(x)f(x),\qquad \qquad Kf(x)=
   \sum_{j=1}^l \langle f,g_j\rangle
g_j(x),
\]
$s(x)$ is a bounded real-valued function and $g_j\in L^2$. Since $K$
has finite rank, Weyl's theorem ensures that $\spec\!_\mathrm{ess}\,
(S+K)=\spec\!_\mathrm{ess}\, S=\spec S= \mathrm{essRange}\,s(x)$. We
aim at finding asymptotic upper bounds for $\delta(n)$ in the limit
$n\to \infty$ in terms of invariants of $s(x)$ and $g_j(x)$, for
suitable $\CL_n$ specified later. The discrete Schr\"odinger
operator studied in \cite{bou1}, \cite[Example I]{lesh} and the
example discussed in \cite[Section~4]{bou2}, all satisfy the
hypothesis of Lemma~\ref{t11} below.

Let $\mathcal{B}_1:=\{e^{inx}\}_{n=-\infty}^\infty$. By declaring
the inner product normalized as $
    \langle f,g\rangle :=\frac{1}{2\pi} \int_{-\pi} ^\pi f(x)
    \overline{g(x)} \mathrm{d}x,
$ clearly  $\mathcal{B}_1$ is an orthonormal basis of $L^2$.
Denote by $\widehat{h}(n):=\langle h,e^{in(\cdot)}\rangle$ the
Fourier coefficient of a complex-valued real function $h$. The
following notation is in place for the examples of this section.
For $q\geq 1$,
\[
   \mathcal{O}^q:=\{h\,:\,{\small
   2\pi-\mathrm{periodic\ in}\ \R}\ \mathrm{s.t.}\
   |\widehat{h}(n)|=O(|n|^{-q})\ \mathrm{as}\ n\to\pm \infty \}.
\]
The well known rules relating Fourier coefficients and
differentiation, easily show that if $h'\in \mathcal{O}^q$, then
$h\in \mathcal{O}^{q+1}$. Hence, $h\in \mathcal{O}^{r+1}$ whenever
$h\in C^r(\R)$ is $2\pi$-periodic. Furthermore, $h\in \mathcal{O}^1$
for functions $h$ of bounded variation, cf. e.g. \cite[\S
2.3.6]{edw}.

\begin{lemma} \label{t11}
Let $M=S+K$ and $\mathcal{L}_n=\mathrm{Span}\,
\{e^{-inx},\ldots,e^{inx}\}$. If $\frac{g_j(x)}{\lambda-s(x)}\in
\mathcal{O}^q$ for all $j=1,\ldots,l$, then we may choose
$\delta(n)=O(n^{-q+1/2})$ as $n\to \infty$ in \eqref{e2}.
\end{lemma}
\proof Since $M=S+K$ is bounded, it is enough estimating the left
hand side of \eqref{e2} for $k=0$. Indeed,
\begin{align*}
  \|\Pi_nM^k\Pi_n\psi-\lambda^k\psi\|&\leq \|\Pi_nM^k\Pi_n\psi-
  \Pi_nM^k\psi\|+\|\lambda^k\Pi_n\psi-\lambda^k\psi\| \\
&\leq (\|M^k\|+|\lambda|^k)\|\Pi_n\psi-\psi\|.
\end{align*}
Now, since $M\psi=\lambda\psi$, then $(\lambda-s(x))\psi(x)=
\sum_{j=1}^l \langle \psi,g_j\rangle g_j(x)$ and so the hypothesis
guarantees that $\psi(x)\in \mathcal{O}^q$. The estimate
\[
\|\Pi_n\psi - \psi\|^2 = \sum_{|j|> n}|\widehat{\psi}(j)|^2 \leq
c_1\sum _{|j|>n} \frac{1}{|j|^{2q}}=O(|n|^{1-2q}), \quad n\to\infty,
\]
completes the proof.

\medskip

In particular we may formulate the following conclusions for $M=S+K$
and $\mathcal{L}_n=\mathrm{Span}\, \{e^{-inx},\ldots,e^{inx}\}$.
\begin{itemize}
\item[a)] If $s(x)$ and $g_j(x)$ are $2\pi$-periodic functions
in $C^r(\R)$, then we may choose $\delta(n)=O(|n|^{-r-1/2})$ in
\eqref{e2}.
\item[b)] If
$s'(x)$ and $g_j'(x)$ are of bounded variation, then we may
choose $\delta(n)=O(|n|^{-3/2})$ in \eqref{e2}.
\item[c)] If $s(x)$ and $g_j(x)$ are of bounded variation,
then we may  choose $\delta(n)=O(|n|^{-1/2})$ in \eqref{e2}.
\end{itemize}
All these assertions follow from the fact that if \linebreak
$\lambda\in \spec\!_\mathrm{disc}\, (S+K)$, then necessarily
$\lambda\not \in \mathrm{essRange}\,s(x)$ and so
$\frac{1}{\lambda-s(x)}$ is of the same degree of smoothness as
$s(x)$. The interesting case in our current discussion is c),
because gaps in the essential spectrum only occur if $s(x)$ is
discontinuous.

\medskip

\subsection{A numerical example} \label{ss2}
Let $M=S+K$ be as in the previous section. Assume that
\begin{equation} \label{e13}
   s(x)=\left\{\begin{array}{lc} -2+\sin (2x),& -\pi<x\leq 0, \\
2+\sin(2x), & 0<x\leq \pi,\end{array} \right.
\end{equation}
$l=1$ and $g_1(x)\equiv \sqrt{2}$, so that $Kf(x)=2\widehat{f}(0)$.
We may find  $\spec M$ in closed form. Clearly
$\spec\!_\mathrm{ess}\,M=[-3,-1]\cup[1,3]$. On the other hand,
according to \cite[Lemma~7]{bou2}, we know that $\lambda$ is an
eigenvalue in the discrete spectrum of $M$ if and only if
\[
   \int_{-\pi}^0 \frac{\mathrm{d} x}{(\lambda+2)-\sin(2x)}+
 \int_0^\pi \frac{\mathrm{d} x}{(\lambda-2)-\sin(2x)} = \pi.
\]
The explicit computation of the integrals reveals that
$\spec\!_\mathrm{disc}\,M$  consists of two eigenvalues, the
variational $\lambda_+\approx 3.5796$ and the pollution-prone
$\lambda_-\approx -0.7674$.

A direct calculation yields
\[
   \widehat{s}(j)=\left\{\begin{array}{lc} \frac{4}{ij\pi},
& j-\mathrm{odd}, \\
(\delta_{2,j}-\delta_{-2,j})\frac{1}{2i}, & j-\mathrm{even},
\end{array} \right.
\]
so that $|\widehat{s}(j)|\sim|j|^{-1}$ as $j\to \pm \infty$ (here
$\delta_{j,l}$ denotes the Kronecker symbol). Thus, by
Lemma~\ref{t11}, we can choose $\delta(n)=O(n^{1/2})$ and, according
to Theorem~\ref{t1}, the existence of $z^{\pm}_n\in\spec P_n$ such
that $|z^{\pm}_n-\lambda_{\pm}|=O(n^{-1/4}) $ in the limit $n\to
\infty$ is predicted.

\begin{table}
\begin{tabular}{|c|c|c|c|c|}
  \hline
 $n$ & $|z^-_n-\lambda_-|$ & $\log(|z^-_n-\lambda_-|)$&
$\log(n)$ & Slope \\
  \hline
190 &0.040879 &  -3.1971 & 5.247 & -0.50849\\
235 &0.036691 &  -3.3052 & 5.4596 & -0.48708\\
280 &0.033689 &  -3.3906 & 5.6348 & -0.50956\\
325 &0.031226 &  -3.4665 & 5.7838 & -0.4876\\
\hline
370 &0.029312 &  -3.5297 & 5.9135 & -0.50963\\
415 &0.027647 &  -3.5882 & 6.0283 & -0.48835\\
460 &0.026291 &  -3.6385 & 6.1312 & -0.50928\\
505 &0.025071 &  -3.6860 & 6.2246 & -0.48918\\
550 &0.024046 &  -3.7278 & 6.3099 & -0.50875\\
\hline
595 &0.023103 &  -3.7678 & 6.3886 & -0.49003\\
640 &0.022292 &  -3.8035 & 6.4615 & -0.50813\\
685 &0.021535 &  -3.8381 & 6.5294 & -0.49086\\
730 &0.020873 &  -3.8693 & 6.593 & -0.50748\\
775 &0.020249 &  -3.8997 & 6.6529 & -0.49167\\
\hline
820 &0.019695 &  -3.9274 & 6.7093 & -0.50682\\
865 &0.019169 &  -3.9545 & 6.7627 & -0.49244\\
 910 &0.018696 & -3.9795 & 6.8134 & -0.50617\\
955 &0.018245 &  -4.0039 & 6.8617 & -0.49318\\
1000& 0.017835 &  -4.0266 & 6.9078 & \\
\hline
\end{tabular}
\caption{\label{ta1} Estimation of the exponent of $\delta(n)$ for
approximations of $\lambda_-$ using $\mathcal{B}_1$.}
\end{table}

From the explicit expression for $\widehat{s}(j)$ it is not
difficult to find $P_n(z)$. In Table \ref{ta1} we report on the
numerical approximations of $|z^-_n-\lambda_-|$ for different values
of $n=190:1000$ and the corresponding pairwise slopes between the
steps $n$ and $n+45$ of the graph $\log|z^-_n-\lambda_-|$ \emph{vs}
$\log (n)$. We have found the numerical data by writing explicitly
$M_n$ and $[M^2]_n$, and computing the eigenvalue  of the companion
matrix of $P_n(z)$ that is nearer to $\lambda_-$. For calculations
we use the standard \texttt{eigs} routine available in the {\scshape
Matlab} package. The last column strongly suggests that
$|z^-_n-\lambda_-|=O(n^{-\alpha})$ for $\alpha \approx 1/2$. Similar
numbers are found for the variational $\lambda_+$.

Table~\ref{ta1} suggests that there is generally a significant gap
between the approximation predicted by \eqref{e20} and the actual
rate of convergence $|z_n- \lambda|\to 0$. An obvious conjecture is
that, perhaps, the bound $\delta(n)^{1/2}$ in Theorem~\ref{t1} can
actually be improved to $\delta(n)$.

\medskip

\subsection{Direct sum of multiplication operators} \label{ss3}
Choosing the right basis is  crucial for achieving efficient
approximation. In the example discussed above, this choice should be
made attending the nature of the symbol, i.e. piecewise continuity.
Below we introduce the correct basis to deal with symbols such as
\eqref{e13}. The following results have obvious extensions to direct
sum of any finite number of operators, i.e. many gaps in the
essential spectrum. In order to keep our notation simple, we only
consider two summing terms.

Assume now that $M=S+K$ is an operator acting on $L^2\oplus L^2$,
where
\begin{equation} \label{e14}
S\begin{pmatrix}f_1(x) \\ f_2(x) \end{pmatrix}:=
\begin{pmatrix} s_1(x) f_1(x)\\s_2(x) f_2(x) \end{pmatrix} ,
\quad K \begin{pmatrix}f_1(x) \\ f_2(x) \end{pmatrix}:=
\sum _{j=1}^l \Big \langle \begin{pmatrix}f_1\\
f_2\end{pmatrix}, G_j \Big\rangle G_j(x),
\end{equation}
$s_1$, $s_2$ are real-valued bounded symbols and
$G_j=(g_{j,1},g_{j,2})^t\in L^2\oplus L^2$. The difference with the
previous model lies on the fact that $\spe _\mathrm{ess}\,
M=\mathrm{essRange}\,s_1(x)\cup \mathrm{essRange}\,s_2(x)$, so gaps
in the essential spectrum may occur even when $s_1$ and $s_2$ are
 continuous and $2\pi$-periodic.

The following lemma is  the natural adaptation of Lemma~\ref{t11} to
the present situation. Let $e_n(x):=( e^{inx} , 0)^t$ and $h_n(x)=(
0, e^{inx})^t$. Then $\{e_n,h_n\}_{n=-\infty}^\infty$ is an
orthonormal basis of $L^2 \oplus L^2$ with the usual inner product.

\begin{lemma} \label{t12}
Let $M=S+K$ where $S$ and $K$ are as in \eqref{e14}. Let
$\CL_n:=\mathrm{Span}\,\{e_{-n},h_{-n},\ldots e_{n},h_{n}\}$. If
$\frac{g_{j,m}(x)}{\lambda-s_m(x)}\in \mathcal{O}^q$ for all
$j=1,\ldots,l$ and $m=1,2$, then we may choose
$\delta(n)=O(n^{-q+1/2})$ as $n\to \infty$ in \eqref{e2}.
\end{lemma}

The proof is analogous to that of Lemma~\ref{t11}.

\medskip

\begin{table}
\begin{tabular}{|c|c|c|c|c|}
  \hline
 $n$ & $|z^-_n-\lambda_-|$ & $\log(|z^-_n-\lambda_-|)$&
$\log(n)$ & Slope \\
  \hline
 12  &   0.037578   &   -3.2813   &    2.4849    &  -2.8383 \\
 18  &   0.011889    &  -4.4322   &    2.8904   &   -3.1543 \\
 24  &  0.0047977    &  -5.3396   &    3.1781   &   -4.9788 \\
 30  &   0.0015796  &    -6.4506   &    3.4012   &   -4.9366 \\
\hline

36 &  0.00064217 &     -7.3507  &     3.5835  &     -7.194  \\
 42 & 0.00021185  &     -8.4596 &      3.7377 &     -6.7381 \\
  48 & 8.6154e-05 & -9.3594  &      3.8712 &      -9.4147 \\
  54 &  2.8424e-05 &    -10.468 &      3.989  &     -8.5396 \\
\hline
 60 &  1.156e-05  &    -11.368  &     4.0943    & -11.635\\
 66 & 3.8138e-06  &    -12.477  &     4.1897    &  -10.341 \\
 72 & 1.5509e-06   &   -13.377  &     4.2767    &   -13.86 \\
 78 & 5.1145e-07    &  -14.486  &     4.3567    &  -12.191 \\
\hline
\end{tabular}
\caption{\label{ta2} Estimation  of $\lambda_-$ using
$\mathcal{B}_2$}
\end{table}

As an application of Lemma~\ref{t12}, we now consider a better
suited basis for estimating the eigenvalues of the symbol discussed
in Section~\ref{ss2}. The crucial observation is that $L^2$ is
isometrically isomorphic to $L^2\oplus L^2$ via the unitary map
\[
   f(x)\stackrel{\mathcal{U}}{\longmapsto}
   \begin{pmatrix} f(\frac{x-\pi}{2}) \\ f(\frac{x+\pi}{2})
   \end{pmatrix}.
\]
The orthonormal basis $\{e_n,h_n\}$ maps under $\mathcal{U}^\ast$
onto an orthonormal basis $\mathcal{B}_2=\{\mathcal{U}^\ast e_n,
\mathcal{U}^\ast h_n\}\subset L^2$ in such manner that
$\mathcal{U}^\ast e_n$ span all functions with support contained in
$[-\pi,0]$ and $\mathcal{U}^\ast h_n$ all those with support in
$[0,\pi]$. When $s(x)$ is as in \eqref{e13},
\[
   \mathcal{U}M\mathcal{U}^\ast\begin{pmatrix}
   f_1(x) \\ f_2(x) \end{pmatrix}=
   \begin{pmatrix} (-2+\sin x)f_1(x) \\ (2+\sin x) f_2(x)
   \end{pmatrix}+2\Big\langle\begin{pmatrix} f_1\\f_2\end{pmatrix}
   ,\begin{pmatrix}1\\1\end{pmatrix}\Big\rangle
   \begin{pmatrix}1\\1\end{pmatrix}.
\]
Then, Lemma~\ref{t12} along with Theorem~\ref{t1} predict that the
basis $\mathcal{B}_2$ would produce sequences approaching to
$\lambda_{\pm}$ whose speed of convergence is faster than any power
of $n$. That is, if $M$ is as in Section~\ref{ss2} and
\[
\CL_n=\mathrm{Span}\,\{\mathcal{U}^\ast e_{-n},\mathcal{U}^\ast
h_{-n},\ldots ,\mathcal{U}^\ast e_{n},\mathcal{U}^\ast h_{n}\},
\]
then for all $q>0$ there is $c(q)>0$ independent of $n$ such that
$|z^\pm_n-\lambda_\pm|\leq c(q)n^{-q}$ for all $n\in \N$.

In Table~\ref{ta2} we report on computations of the nearest point in
$\spec P_n$ to $\lambda_-$ and list the values analogous to
Table~\ref{ta1}. In this case we only consider $n=12:78$. Notice
that the approximation of $\lambda_-$ using $\mathcal{B}_2$ at
$n=12$ is already more accurate than the one obtained using the
basis $\mathcal{B}_1$ at $n=1000$. The super-polynomial speed of
approximation predicted by Lemma~\ref{t12} and Theorem~\ref{t1} is
already evidenced by the first few entries of the last column.

\subsection{Schr\"odinger operators with band gap essential spectrum}
\label{ss4}

For our last example we consider an unbounded operator. Let $M$ be
the one-dimensional Schr\"odinger operator
\begin{equation} \label{e22}
   Mf(x)=-\partial_x^2f(x)+V(x)f(x), \qquad \qquad f\in \dom M:=W^{2,2}(\R)
\end{equation}
acting on $L^2(\R)$, where the potential $V=V_1+V_2$, $V_1\in
W^{2,2}(\R)$ and $V_2\in W^{2,\infty}(\R)$. Here $W^{n,p}(\R)$
denotes the $n^\mathrm{th}$ derivative $L^p$ Sobolev space: $V\in
W^{n,p}(\R)$ if and only if $V\in L^p(\R)$ and $\partial_x^q V\in
L^p(\R)$ for $q=1,\ldots,n$.

Since, in particular, $V\in L^2(\R)+L^\infty(\R)$, multiplication by
$V$ is $\partial_x^2$ bounded with relative bound $0$
\cite[Theorem~XIII.96]{resi}. Recall that the operator $L$ is said
to be $L_0$ bounded, if and only if $\dom L_0\subseteq \dom L$ and
\[
   \|Lf\|\leq \alpha \|L_0f\|+\beta \|f\|, \qquad f\in \dom L_0
\]
for suitable non-negative constants $\alpha$ and $\beta$.  The
infimum of all $\alpha$ allowing the above for some $\beta$ is
called the relative bound of $L$ (with respect to $L_0$). By virtue
of Kato-Rellich theorem \cite[Theorem~X.12]{resiv2}, the above
choice of $\dom M$ guarantees $M=M^\ast$. Notice that $\dom
M^2=W^{4,2}(\R)$.

Furthermore, since $V_1$ is relatively compact with respect to
$-\partial_x^2$ \cite[Theorem~XIII.15]{resi} and $V_2$ is bounded,
$V_1$ is also relatively compact with respect to
$-\partial_x^2+V_2$. Hence, by Weyl's theorem, the essential
 spectrum of $H$ is characterized completely by $V_2$, that
is $\spe_{\mathrm{ess}} M= \spe_\ess(-\partial_x^2+V_2)$. If $V_2$
is periodic  \cite[Chapter XIII.16]{resi}, the essential spectrum is
bounded below, it lies in bands and it extends to infinity. In
general, $V_1$ may produce non-empty discrete spectrum.

 We put
\begin{equation} \label{e24}
   \CL_n=\mathrm{span} \{\phi_0,\ldots,\phi_n\}
\end{equation}
where $\phi_j(x)=c_kh_j(x)e^{-x^2/2}$, $c_j=\frac{1}{\sqrt{2^j
j! \sqrt{\pi}}}$, and $h_j(x)$ is the $j^{\mathrm{th}}$ Hermite  polynomial
given by Rodriguez's formula
\[
   h_j(x)=(-1)^j e^{x^2}(\partial_x^j e^{-x^2}).
\]
The choice of $c_j$ ensures that $\|\phi_j\|=1$ for all $j=0,1,\ldots$.
It is well known
that $\phi_j(x)$ are the eigenvectors of the
quantum mechanical harmonic oscillator $H:=-\partial_x^2+x^2$, so that
$\{\phi_j\}_{j=0}^\infty$ is an orthonormal basis of $L^2(\R)$.

\begin{theorem} \label{t21}
Let $M$ be given by \eqref{e22} and $\CL_n$ be given by \eqref{e24}.
Then for all $\lambda\in \spe_{\mathrm{disc}} M$ there is
$z_n\in\spec P_n$ such that $z_n\to \lambda$ as $n\to\infty$.
Furthermore, if $\partial^q_xV(x)$ is continuous and bounded for
some $q\geq 3$, then $|z_n-\lambda|=o(n^{\frac{2-q}{4}})$ as $n\to
\infty$.
\end{theorem}

The proof of this result will be a consequence of
three technical lemmas.

\begin{lemma} \label{t20}
If $\partial^q_xV(x)$ is continuous and
bounded for
some $q\geq 1$, then $\psi\in C^{q+2}(\R)$
and $|\partial_x^{(q+2)}\psi(x)|<d_1e^{-a_1|x|}$, $x\in \R$,
for suitable constants
$d_1,a_1>0$.
\end{lemma}
\proof Since $V\in W^{2,2}(\R)$, then $V'(x)$ is a continuous
function. By hypothesis, $\psi''=(V-\lambda)\psi$, showing that
$\psi''(x)$ is also continuous. By repeating recursively
this argument $q$ times, the equality
\begin{equation} \label{e21}
 \partial^{(q+2)}_x\psi = \psi\partial_x^q(V-\lambda)+\ldots+(V-\lambda)
\partial_x^q\psi
\end{equation}
implies the existence and the needed continuity of
$\partial_x^{(q+2)}\psi(x)$.

By virtue of the results of \cite[\S C.3]{250}, the hypothesis on
$V$ ensures that
\begin{equation} \label{e23}
   |\psi(x)|< d_2e^{-a_2|x|}, \qquad x\in \R,
\end{equation}
for some positive constants $d_2,a_2$. Since $V(x)-\lambda$
is continuous and bounded, then
$|\psi''(x)| < \tilde{d}_2e^{-\tilde{a}_2|x|}$.
By repeating  recursively this argument $q$ times,
\eqref{e21} and the hypothesis, yield
the desired upper bound for $|\partial_x^{(q+2)}\psi(x)|$.

\medskip

\begin{lemma} \label{t22}
Let $M$ be given by \eqref{e22}. Then any eigenfunction of $M$ lies
in $\mathrm{Dom}\, H^2$.
\end{lemma}
\proof From the inclusion $\psi\in \dom M^2=W^{4,2}(\R)$ and
\eqref{e23} (which do not require any smoothness property on $V$),
it follows that $\psi\in W^{2,2}(\R)\cap \dom(x^2)=\dom H$ and
$(-\psi''+x^2\psi)\in \dom H$. This ensures the desired property.
The only non-trivial facts in the latter assertions are, perhaps,
the inclusions $x^2\psi\in W^{2,2}(\R)$ and $\psi''\in \dom(x^2)$.
The first follows from the second, by differentiating twice the term
$x^2\psi$ and noticing that $(x\psi')'\in L^2(\R)$. The second one
is achieved as follows. Since $M\psi$ is an eigenvector of $M$, $
   |-\psi''(x)+V(x)\psi(x)|\leq d_3 e^{-a_3|x|}.
$
Then
\begin{align*}
  \|x^2 \psi''\|& \leq \|x^2(-\psi''+V\psi)\|+\|x^2V\psi\| \\
   & \leq d_4+ \left(\int x^4 |V(x)|^2\,|\psi(x)|^2
   \mathrm{d} x\right)^{1/2} \\
   & \leq d_4 + d_2^2\left( \int x^4 e^{-2a_2|x|} |V(x)|^2 \mathrm{d}
x\right)^{1/2}.
\end{align*}
The latter integral is bounded because of $V\in L^2(\R)+L^\infty(\R)$,
therefore \linebreak $\psi''\in \dom(x^2)$.

\medskip

\begin{lemma} \label{t23}
Let $M$ be given by \eqref{e22}. Then both $M$ and $M^2$ are $H^2$
bounded.
\end{lemma}
\proof Multiplication by $V$ is $\partial^2_x$ bounded with bound 0.
Then $M$ is $\partial_x^4$ bounded with relative bound 0. Similarly,
from the identity
\[
   (-\partial_x^2+V)^2u =
   \partial_x^4u-2V\partial_x^2u-2V'\partial_xu+(V^2-V'')u
\] and the fact that $V,V',V''$ lie in $L^2(\R)+L^\infty(\R)$, one may
deduce that $M^2$ is $\partial_x^4$ bounded with relative bound 1.
Then the proof reduces to showing that $\partial_x^4$ is $H^2$
bounded. For this, let
\[
    A:=2^{-1/2}(x+\partial_x)\qquad \mathrm{and} \qquad
    A^\ast = 2^{-1/2}(x-\partial_x).
\]
Then $(-\partial_x^2+x^2)=2(AA^\ast-1)$ and
$\partial_x^4=(A-A^\ast)^4$. Thus the desired property follows
from the identity (cf.\cite[eq.(X.28)]{resiv2})
\[
    \|A^{\#_1}\cdots A^{\#_q}u\|\leq c
    \|(-\partial_x^2+x^2)^{q/2}u\|, \qquad q=1,2,\ldots
\]
where $A^{\#_k}$ is either $A$ or $A^\ast$. This is
easily shown by induction and using
the estimate
\[
  \|(-\partial_x^2+x^2)^{p/2}u\|\leq\|(-\partial_x^2+x^2)^{q/2}u\|,
  \qquad p<q.
\]
This completes the proof of the lemma.

\medskip

We may now complete the proof of Theorem~\ref{t21}. Let $\phi:=H^2
\psi$. The crucial point lies in the following observation.
According to Lemma~\ref{t23}, there exist non-negative constants
$\alpha$ and $\beta$, such that
\begin{align*}
\|\Pi_nM^k\Pi_n\psi-\lambda^k\psi\|& \leq \|\Pi_nM^k\Pi_n\psi-
\Pi_nM^k\psi\|+\|\Pi_nM^k\psi-\lambda^k\psi\| \\
&\leq \|M^k(\Pi_n\psi-\psi)\|+\|\lambda^k(\Pi_n\psi-\psi)\| \\
&\leq \alpha
\|H^2(\Pi_n\psi-\psi)\|+(\beta+|\lambda|^k)\|\Pi_n\psi-\psi\|\\&=
\alpha\|\Pi_n\phi-\phi\|+(\beta+|\lambda|^k)\|\Pi_n\psi-\psi\|.
\end{align*}
This ensures directly the first part of the theorem. For the second
part we only require estimating $\|\Pi_n\psi-\psi\|$ and
$\|\Pi_n\phi-\phi\|$ as $n\to \infty$.

A straightforward application of lemma~\ref{t20} shows that, if
\linebreak $\partial^q_xV(x)$ is continuous and bounded for $q\geq
2$, then for any set of \linebreak polynomials
$\{p_0(x),\ldots,p_{q+2}(x)\}$, the functions $\sum_{j=0}^{q+2}
p_j(x) \psi^{(j)}(x)$ and \linebreak $\sum_{j=0}^{q-2} p_j(x)
\phi^{(j)}(x)$  are continuous and square integrable. By virtue of
the fundamental relation
\begin{equation*}
    h_j(x)=\frac{h'_{j+1}(x)}{2(j+1)}, \qquad \qquad j=0,1,\ldots,
\end{equation*}
integration by parts yields
\begin{align*}
\langle \psi,\phi_k\rangle & = c_k\int_{-\infty}^\infty \psi(x)
h_k(x) e^{-x^2/2} \ud x = c_k\int_{-\infty}^\infty \psi(x)
\left(\frac{h'_{k+1}(x)}{2(k+1)}\right) e^{-x^2/2} \ud x \\
& = -\frac{c_k}{2(k+1)} \int_{-\infty}^\infty [\psi(x)
e^{-x^2/2}]' h_{k+1}(x) \ud x \\
& = -\frac{c_k}{2(k+1)} \int_{-\infty}^\infty [\psi'(x)
-x\psi(x)] e^{-x^2/2}h_{k+1}(x) \ud x = \ldots \\
& = \frac{(-1)^{q+2}c_k}{2^{q+2}(k+1)\cdots(k+q+2)}
\int_{-\infty}^\infty \Big[\sum_{j=0}^{q+2}p_j(x)\psi^{(j)}(x)\Big]
e^{-x^2/2}h_{k+q+2}(x) \ud x \\
& = \frac{(-1)^{q+2}\langle \zeta,\phi_{k+q+2}\rangle
}{2^{\frac{q+2}{2}}\sqrt{(k+1)\cdots(k+q+2)}},
\end{align*}
where $\zeta(x)=\sum_{j=0}^{q+2}p_j(x)\psi^{(j)}(x)$ is continuous
and square integrable. Thus
\begin{align*}
\|\Pi_n\psi-\psi\|^2 &= \sum_{k=n+1}^\infty |\langle \psi,\phi_k\rangle|^2
\\ & \leq  \sum_{k=n}^\infty \frac{|\langle \zeta,\phi_{k+q+2}\rangle|^2}
{2^{q+2}(k+1)\cdots(k+q+2)} \\ &\leq
\frac{1}{2^{q+2}n^{q+2}}\sum |\langle \zeta,\phi_{k+q+2}\rangle|^2
= o(n^{-(q+2)})
\end{align*}
as $n\to \infty$. A similar calculation with $\phi$ instead of
$\psi$, gives
\[
    \|\Pi_n\phi -\phi\|^2= o(n^{-(q-2)}), \qquad n\to \infty,
\]
completing the proof of Theorem~\ref{t21}.

\medskip

We remark that Theorem~\ref{t21} may be extended to higher dimension
without much effort.

\section{Distance to singularity, the spectral function}

In order to simplify the notation in the proof of Theorems~\ref{t1}
and \ref{t5}, we first discusses the notion of distance from $P(z)$
to the nearest singular matrix. For convenience, all the ideas in
this section refer to general matrix polynomials of arbitrary degree
$P(z)=\sum_{k=0}^{m} A_kz^k$, where $A_k\in \C^{j\times j}$ and
$\det A_m\not= 0$. Some of the results below are related to the
recent electronic manuscript by Davies \cite{mao} and the paper by
Lancaster and Psarrakos \cite{laps}.

For all $z\in \C$, we define
\[
   \sigma_P(z):= \inf_{v\not = 0} \frac{\|P(z)v\|}{\|v\|}.
\]
Four different characterizations of this function are,
\begin{itemize}
\item[a)] $\sigma_P(z)$ equals the smallest singular value of
$P(z)$,
\item[b)] $\sigma_P(z)^{-1}=\|P(z)^{-1}\|$ whenever $\det
P(z)\not=0$,
\item[c)] $\sigma_P(z)^{-1}=\sup_{u,v\not=0}\frac{\Re\langle
P(z)^{-1}u,v\rangle}{\|u\|\|v\|}$ whenever $\det
P(z)\not=0$,
\item[d)] $\begin{array}{rl} \sigma_P(z)&=\min \{\|E\|: \det[P(z)+E]=0,
E\in \C^{j\times j}\}
\\ &= \min \{\|E\|: \det[P(z)+E]=0, \mathrm{Rank}\, E=1\}. \end{array}$
\end{itemize}
Since
\[
   \spec P=\{z\in \C:\sigma_P(z)=0\},
\]
we may call the scalar quantity $\sigma_P(z)$ a ``spectral
function'' of the matrix polynomial $P$.

The proofs of a), b) and c) are straightforward and property d)
holds trivially when $\sigma_P(z)=0$. In order to prove d) when
$\sigma_P(z)\not=0$, we use b).  If $\|E\|<\|P(z)^{-1}\|^{-1}$, then
$\|P(z)^{-1}E\|<1$ and $P(z)+E=P(z)(I+P(z)^{-1}E)$ so that
$\det(P(z)+E)\not=0$. This shows that $\|P(z)^{-1}\|^{-1}$ can not
be greater than the right sides of d). Conversely, let $v\in \C^j$
be such that $\|P(z)^{-1}v\|=\|P(z)^{-1}\|$ and $\|v\|=1$, and let
$u:=\|P(z)^{-1}\|^{-1}P(z)^{-1}v$. Then $\|u\|=1$. Put
$Ew:=-\|P(z)^{-1}\|^{-1}\langle w,u\rangle v$ for all $w\in \C^j$.
Then the linear operator $E$ has rank one, it satisfies
$\|E\|=\|P(z)^{-1}\|^{-1}$ and $Eu=-\|P(z)^{-1}\|^{-1}v$. Thus
$(P(z)+E)u=0$ so that $\det[P(z)+E]=0$. This ensures that the bottom
right side of d) is not greater than $\|P(z)^{-1}\|^{-1}$.

\medskip

The following lemma will play a fundamental role in the sequel. See
\cite[Lemma~2]{bou2} and \cite[Section~3.4]{mao}.

\begin{lemma} \label{t2}
The non-negative function $\sigma_P(z)$ is Lipschitz continuous in
every compact subset of the complex plane. Furthermore
$\sigma_P(z)^{-1}$ is subharmonic in $\C\setminus \spec P$, so that
$z_0$ is a local minimum of $\sigma_P$ if and only if $z_0\in \spec
P$.
\end{lemma}
\proof The first part of the lemma follows easily from the triangle
inequality once we have proven the estimate
\begin{equation} \label{e4}
|\sigma_P(z)-\sigma_P(w)|\leq \|P(z)-P(w)\|
\qquad \mathrm{for\ all}\ z,w\in \C.
\end{equation}
We show the latter by considering two separated cases. If $z\in
\spec P$ and $w\not \in \spec P$, since $P(z)=P(w)+(P(z)-P(w))$ is
not invertible, d) ensures \eqref{e4}. If $z,w\not \in \spec P$,
then \[P(z)^{-1}-P(w)^{-1}= P(z)^{-1}[P(w)-P(z)]P(w)^{-1}.\] By the
triangle inequality,
\[
    \left| \|P(z)^{-1}\|-\|P(w)^{-1}\|\right|
    \leq\|P(z)^{-1}\|
    \|P(w)-P(z)\|\|P(w)^{-1}\|
\]
and so \eqref{e4} is consequence of b). The second part of the lemma
follows from c) and the elementary properties of subharmonic
functions.

\medskip

We now consider weighted perturbations of $P(z)$ in the sense
studied recently in \cite{hiti} and \cite{laps}. We will  require
these  results in the proof of Theorem~\ref{t5}. Below,
$Q(z)=\sum_{k=0}^{m}E_k z^k$ shall always refer to a small
perturbation of $P$. Let $\varepsilon\geq 0$ and
$\bar{w}=(w_0,\ldots,w_n)$ where $w_k\geq 0$. The (weighted)
$\varepsilon$-pseudospectra of $P$ are the set given by {\small
\[
   \spec\!_{\varepsilon,\bar{w}}\, P:=
   \Big\{z\in \C: \begin{array}{l}
   \det [P(z)+Q(z)]=0, \,\mathrm{for\ some}\ Q(z)=\sum_{0}^{m}E_k z^k
   \\ \mathrm{such\ that}\ \|E_k\|\leq w_k \varepsilon,\,
   k=0,\ldots,m \end{array}\Big\}.
\]}
This definition was studied by Higham and Tisseur in \cite{hiti} and
it extends the standard definition of pseudospectra of a matrix $A$
in the obvious manner, by considering $P(z)=(z-A)$ and
$\bar{w}=(1,0)$. The weight $\bar{w}$ is introduced in order to
allow freedom in controlling the perturbation of each individual
coefficient of $P$, e.g. this may be given in an absolute sense
($w_k=1$) or in a relative sense ($w_k=\|A_k\|$). Observe that
$\spec\!_{0,\bar{w}} P=\spec P$.

By combining the remarkable result \cite[Lemma~2.1]{hiti} with b),
we achieve the useful characterization
\[
   \spec\!_{\varepsilon,\overline{w}}\, P=\{z\in \C:\sigma_P(z)\leq
   \varepsilon (w_0+w_1|z|+\ldots+w_m|z|^{m}) \}.
\]
Notice that we do not require $\|E_m\|$  to be small, so
$\spec\!_{\varepsilon,\bar{w}}\,P$ is not guaranteed to be bounded.
This shall not be an important point here. Necessary and sufficient
conditions for these set to be bounded may be found in
\cite[Theorem~2.2]{laps}.

The following lemma is relevant in the proof of Theorem~\ref{t5}.
This has also been established by Lancaster and Psarrakos in the
more sophisticated \cite[Theorem~2.3]{laps}.

\begin{lemma} \label{t4}
Let $\varepsilon>0$. Let $\Omega$ be a connected component of
$\spec\!_{\varepsilon,\bar{w}}\,P$ such that $\Omega\cap \spec
P\not=\varnothing$. If $\|E_k\|\leq w_k \varepsilon$ for all
$k=1,\ldots,m$, then $\Omega \cap \spec (P+Q) \not= \varnothing$.
Furthermore $P+Q$ has the same number of eigenvalues in $\Omega$,
counting multiplicity, than $P$ has.
\end{lemma}
\proof Let $P_\delta(z):=\sum_{k=0}^m(A_k+\delta E_k)z^k$ for $0\leq
\delta\leq 1$. Then $P_0(z)=P(z)$ and $P_1(z)=P(z)+Q(z)$. Since
these eigenvalues are the zeros of a certain family of scalar
polynomials whose coefficients depend continuously upon $\delta$,
the eigenvalues of $P_\delta(z)$ depend continuously on $\delta$.
Let $\mu_\delta\in \spec P_\delta$, be such that $\mu_0\in \Omega$.
Since $\|\delta E_k\|\leq \|E_k\|\leq w_k \varepsilon$, then
$\mu_\delta \in \spec\!_{\varepsilon,\bar{w}}\,P$ for all $0\leq
\delta \leq 1$. Being a continuous trajectory, necessarily
$\mu_\delta \in \Omega$ for all $0\leq \delta \leq 1$.


\section{Proof of theorems~\ref{t1} and \ref{t5}}

Throughout this section we shall write $\sigma_n\equiv
\sigma_{P_n}$.
\medskip

Let $b_1:=2+2|\lambda|+3|\lambda|^2$. By virtue of the triangle
inequality,
\begin{align*}
\|\Pi_n&(\lambda-M)^2\Pi_n\psi\| = \|\lambda^2\Pi_n\psi-2\lambda
\Pi_nM\Pi_n\psi+\Pi_nM^2\Pi_n\psi\| \\ & \leq
\|\Pi_nM^2\Pi_n\psi-\lambda^2\Pi_n\psi\|+\|2\lambda^2\Pi_n\psi-2\lambda
\Pi_nM\Pi_n \psi\|
\\ &\leq \|\Pi_nM^2\Pi_n\psi-\lambda^2\psi\|+
2|\lambda|\,\|\lambda\psi-\Pi_nM\Pi_n\psi\|+3|\lambda|^2
\|\psi-\Pi_n\psi\|\\ &< (b_1-1/2) \delta(n).
\end{align*}
Then, from the definition of $\sigma_n$, by choosing $v=\Pi_n\psi$,
we achieve the following.

\begin{lemma} \label{t7}
There exists $N_1>0$ such that
$
  \sigma_n(\lambda)< b_1 \delta(n)
$
for all $n>N_1$.
\end{lemma}

\medskip

The proof of both theorems depends heavily on applying the second
part of Lemma~\ref{t2}. For this we combine Lemma~\ref{t7} with
lower bounds for $\sigma_n(z)$ along the boundary of small
neighbourhoods of $\lambda$. Our first step towards computing these
bounds, involves finding estimates for $\|P_n(z)^{-1}\|$ in
subspaces that are orthogonal to the eigenspace associated to
$\lambda$.

Let $\CE:=\mathrm{span}\,\{ \psi \in \dom M : M\psi=\lambda \psi\}$,
$d_\lambda:=\dim \CE$ and $\Pi_\CE$ be the orthogonal projection
onto $\CE$. Let $\tilde{\mu}=\pm\mu$ be such that
$\lambda-\tilde{\mu}\in \spec M$ and
\[
    A:= M(I-\Pi_\CE)+(\lambda-\tilde{\mu}) \Pi_\CE, \qquad \dom
A=\dom M.
\]
Then, by construction,
$A=A^\ast$ and
$
   \spec A= \spec M \setminus \{\lambda\}.
$ For $z\in \C$, we write $A(z):=(z-A)^2$ with $\dom A(z)=\dom M^2$,
$K(z)=\tilde{\mu}(2\lambda-2z-\tilde{\mu})\Pi_\CE$ so that
\[
   A(z)x=(z-M)^2x-K(z)x\qquad \mathrm{for} \qquad x\in \dom M^2,
\]
$A_n(z):=\Pi_nA(z)|\CL_n$, $K_n(z)=\Pi_nK(z)|\CL_n$ and
$A_n(z)^{-1}=[A_n(z)]^{-1}$.

Let $D:=\{|z-\lambda|\leq \mu/4\}$. Then $D$ does not intersect
$\spec A$ and $\lambda$ is the closest point in $\spec M$ to the
boundary of $D$. Furthermore, if $z\in D$, then
\begin{align*}
    \inf_{w\in \R \setminus [\lambda-\mu,\lambda+\mu]}
    \mathrm{Re}\,(z-w)^2& =\inf_{w\in \R \setminus
    [\lambda-\mu,\lambda+\mu]}
     (w-\lambda-\Re z)^2-(\Im z)^2 \\
     & \qquad \qquad \geq \mu^2/2>0.
\end{align*}
Put $\gamma:=\mu^2/2$. Since $\langle \Pi_n(z-A)^2|\CL_n v,v\rangle
=\langle (z-A)^2 v,v\rangle$ for all $v\in \CL_n $ and $A(z)$ is a
normal operator, then the following inclusion for the numerical
range of $A_n(z)$ holds for all $z\in D$,
\begin{align*}
  \num A_n(z) &\subseteq \num A(z) \subseteq \overline{\mathrm{Conv}
  [\spec (z-A)^2]} \\
  &\subseteq \{(z-w)^2:w\in[\lambda- \mu,\lambda+\mu]\}
  \subseteq \{ \Re (z)\geq \gamma>0\}.
\end{align*}
Thus $A_n(z)$ is invertible and
\[
   \|A_n^{-1}(z)\|\leq [\dist (0,\num A_n(z))]^{-1}\leq
   \gamma^{-1}
\]
for all $z\in D$ and $n\in \N$.

Let
\[
   b_2:=1+2(\mu/4+|\lambda|)+(\mu/4+|\lambda|)^2
   +2\mu(4|\lambda|+3\mu/2).
\]
The triangle inequality ensures that
\begin{align*}
  \|A_n(z)&\Pi_n\psi-A(z)\psi\|\\& =\|\Pi_n(z-M)^2\Pi_n\psi-(z-M)^2\psi
  -\mu(2\lambda-2z-\mu)[\psi-\Pi_n\Pi_\CE \Pi_n\psi]\| \\
  &\leq \|\Pi_nM^2\Pi_n\psi-\lambda^2\psi\|+2|z|\|\Pi_nM\Pi_n\psi-\lambda
  \psi\|+\\& \qquad
  \big(2|\mu(2\lambda-2z-\mu)|+|z|^2\big)\|\psi-\Pi_n\psi\| \\
  &\leq (1+2|z|+|z|^2+2\mu|2\lambda-2z-\mu|)\delta(n)\,\leq
  \,b_2\delta(n)
\end{align*}
for all $z\in D$ and $n\in \N$. Let
$\nu:=(\lambda-\tilde{\mu}-z)^2$, so that $A(z)\psi=\nu\psi$. Then
\begin{align*}
\|A_n(z)^{-1}\Pi_n\psi&-A(z)^{-1}\psi\|=\|A_n(z)^{-1}\Pi_n\psi-\nu^{-1}\psi\|
\\ &\leq \|A_n(z)^{-1}\Pi_n\psi-\nu^{-1}\Pi_n\psi\|+
     \|\nu^{-1}\Pi_n\psi-\nu^{-1}\psi\| \\
     & \leq \gamma^{-1}|\nu|^{-1}\|\nu \Pi_n\psi-A_n(z)\Pi_n\psi\|+
     |\nu|^{-1} \|\Pi_n\psi-\psi\| \\ &= \gamma^{-1}|\nu|^{-1}
   \|\Pi_nA(z)\psi-A_n(z)\Pi_n\psi\|+|\nu|^{-1} \|\Pi_n\psi-\psi\| \\
   & \leq \gamma^{-1}|\nu|^{-1}
   \|A(z)\psi-A_n(z)\Pi_n\psi\|+|\nu|^{-1} \|\Pi_n\psi-\psi\| \\
   &\leq (\gamma^{-1}b_2+1)|\nu|^{-1}\delta(n) \\
   &\leq (4(\gamma^{-1}b_2+1)\mu/3)\delta(n)=:b_3\delta(n)
\end{align*}
for all $z\in D$ and $n\in \N$. Hence, a straightforward computation
shows that
\begin{equation} \label{e8}
   \|A_n(z)^{-1}\Pi_nK(z)-A(z)^{-1}K(z)\|<d_\lambda b_3 \delta(n)
\end{equation}
for all $z\in D$ and $n\in \N$.

The following estimate, which is a direct consequence of the
definition of $K(z)$, will be needed below: \[\|K(z)\|\leq
\mu(2|z-\lambda|+\mu)\leq 3\mu^2/2=:b_4 \qquad \mathrm{for\ all}\
z\in D.\]

\medskip


\subsection{Proof of Theorem~\ref{t1}}
By virtue of Lemma~\ref{t2}, the only local minima of $\sigma_n(z)$
are those points in $\spec P_n$. Then the proof of Theorem~\ref{t1}
reduces to finding $\tilde{b}>0$ and $N_2>0$, both independent of
$z$ and $n$, such that
\begin{equation} \label{e25}
\sigma_n(z)>\sigma_n(\lambda) \qquad \mathrm{for\ all} \quad
z\in\{|z-\lambda|=\tilde{b} \delta(n)^{1/2}\} \quad \mathrm{and}
\quad n>N_2.
\end{equation}
Since $\spec P_n$ is finite, the existence of $b>0$ as we require in
the thesis part of the theorem becomes obvious.

Inequality \eqref{e25} is a consequence of Lemma~\ref{t7} and the
following.

\begin{lemma} \label{t6}
let $\tilde{b}:=(4b_4\max\left\{d_\lambda b_3,\gamma^{-1}b_1\right\}
+1)^{1/2}$ and let $N_2\in \N$ be such that
\begin{equation*}
    \delta(n)< \min \{\tilde{b}^{-2}(\mu/4)^2,
    (4d_\lambda b_3)^{-1}, \gamma(4b_1)^{-1}\}, \qquad n>N_2.
\end{equation*}
Then $\sigma_n(z)>b_1 \delta(n)$ for all
$|z-\lambda|=\tilde{b}\delta(n)^{1/2}$ and $n>N_2$.
\end{lemma}
\proof Throughout the proof we assume that $n>N_2$ and
$|z-\lambda|=\tilde{b}\delta(n)^{1/2}$. Since
$\tilde{b}\delta(n)^{1/2}<\mu/4$, then $A(z)+K(z)=(z-M)^2$ is
invertible and
\begin{equation} \label{e3}
\begin{aligned}
\|[I+&A(z)^{-1}K(z)]^{-1}\| = \|(z-M)^{-2}A(z)\|\\ & =\|(z-M)^{-2}
[(z-M)^2-K(z)]\| \\ &=\|I-(z-M)^{-2}K(z)\| \leq 1+
b_4\|(z-M)^{-2}\|\\ &\leq 1+b_4[\dist (z,\spec M)]^{-2}
=\frac{\tilde{b}^2\delta(n)+b_4}{\tilde{b}^2\delta(n)}
=:[c(n)]^{-1}.
\end{aligned}
\end{equation}
Thus, given $v\in \CL_n$,
\begin{align*}
  c&(n)\|v\| \leq
   \|v+A(z)^{-1}K(z)v\| \\
   & \leq \|v+A_n(z)^{-1}\Pi_nK(z)v\| +
   \|A_n(z)^{-1}\Pi_nK(z)-A(z)^{-1}K(z)\|\|v\|.
\end{align*}
Since $\tilde{b}^2> 4b_4d_\lambda b_3$ and $\delta(n)<(4d_\lambda
b_3)^{-1}$, then $
   \tilde{b}^2>\frac{2b_4d_\lambda b_3}{1-2\delta(n)d_\lambda
   b_3}
$ so that \eqref{e8} and an easy calculation yield $
   \|A_n(z)^{-1}\Pi_nK(z)-A(z)^{-1}K(z)\|<c(n)/2.
$
Hence
\begin{align*}
  \frac{c(n)}{2}\|v\| & \leq \|v+A_n(z)^{-1}\Pi_nK(z)v\| \\
    & \leq \|A_n(z)^{-1}\| \|A_n(z)v+\Pi_nK(z)v\| \\
    &\leq \gamma^{-1} \|\Pi_n(z-M)^2|\CL_nv\|.
\end{align*}
Therefore $\Pi_n(z-M)^2|\CL_n$ is invertible and
\[
  \sigma_n(z)=\|(\Pi_n(z-M)^2|\CL_n)^{-1}\|^{-1} \geq
  \frac{\gamma c(n)}{2}.
\]
Since $\tilde{b}^2>\gamma^{-1}4b_4b_1$ and $\delta(n)<\gamma
(4b_1)^{-1}$, it is easy to see that
$\tilde{b}^2>\frac{2b_4b_1}{\gamma-2b_1\delta(n)}. $ Thus $
   \sigma_n(z)>b_1\delta(n)
$
as claimed in the lemma.

\medskip

This completes the proof of Theorem~\ref{t1}. Explicit expressions
for $\tilde{b}$ in terms of $\lambda$ and $\mu$ might be useful in
applications. A direct substitution yields
$\tilde{b}=(\max\left\{b_5,12b_1\right\}+1)^{1/2} $ where
\[
b_5:=d_\lambda \frac{16 + 16|\lambda|^2 + 8\mu + 57\mu^2 +
8|\lambda|(4 + 17\mu)}{\mu}.
\]


\subsection{Proof of Theorem~\ref{t5}} The proof will be
a consequence of lemmas~\ref{t4}, \ref{t7} and the following.

\begin{lemma} \label{t8}
For all $0<\delta<\mu/4$, there exists $N_3>0$ independent of $z$ or
$n$, such that
\[
    \sigma_n(z)>\frac{\delta^2 \gamma}{2(\delta^2+b_4)}
\]
for all $n>N_3$ and $\delta\leq |z-\lambda|\leq \mu/4$.
\end{lemma}
\proof Let $C:=\{\delta\leq |z-\lambda|\leq \mu/4\}$. Throughout the
proof we assume $z\in C$. Since  $C\subset D$, then
$A(z)+K(z)=(z-M)^2$ is invertible. By virtue of the third line of
\eqref{e3},
\begin{equation*}
\|[I+A(z)^{-1}K(z)]^{-1}\| \leq 1+b_4[\dist (z,\spec M)]^{-2} \leq
\frac{\delta^2+b_4}{\delta^{2}}=:c_2^{-1}.
\end{equation*}
Let $N_3>0$ be such that $d_\lambda b_3\delta(n)<c_2/2$ for all
$n>N_3$. Then a straightforward computation along with \eqref{e8}
yield
\[
   \|A_n(z)^{-1}\Pi_nK(z)-A(z)^{-1}K(z)\|<\frac{c_2}{2}
\]
for all $n>N_3$. Given $v\in \CL_n$, then
\begin{equation*}
  c_2\|v\| \leq \|v+A_n(z)^{-1}\Pi_nK(z)v\| +
   \|A_n(z)^{-1}\Pi_nK(z)-A(z)^{-1}K(z)\|\|v\|.
\end{equation*}
Thus $c_2\|v\|/2\leq \|P_n(z)v\|\gamma^{-1}$, so that $P_n(z)$ is
invertible and $\sigma_n(z)\geq \gamma c_2/2$ for all $n>N_3$.

\medskip

\begin{figure}[t]
\epsfig{file=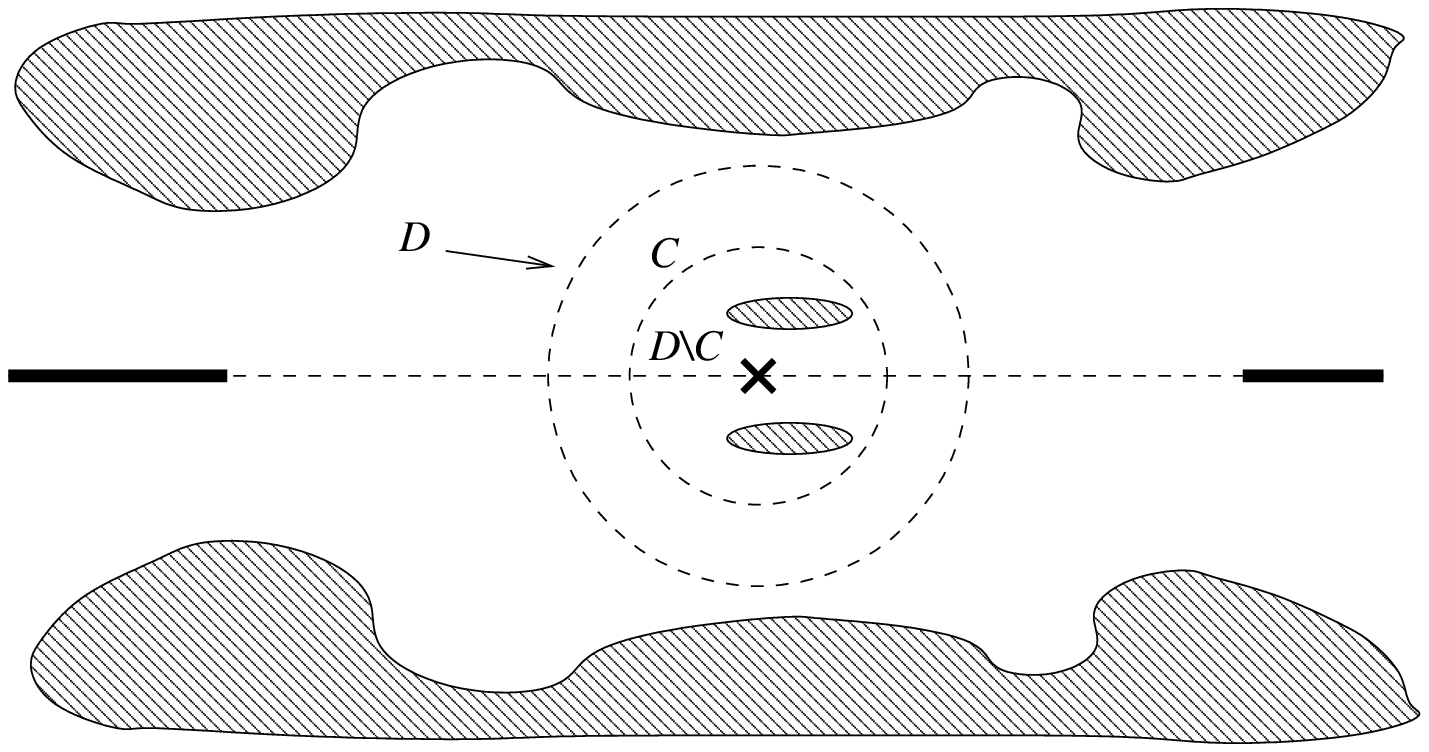, height=2.3in} \caption{Proof of
Theorem~\ref{t5}. The shaded region represents
$\spe_{\varepsilon,\tilde{w}} P_n$ which is symmetric with respect
to the real axis (dotted line). The dark intervals represent
$\spe_{\mathrm{ess}}M$ and the cross the eigenvalue $\lambda$.}
\end{figure}

We now complete the proof of Theorem~\ref{t5}. Notice that the case
$\varepsilon=0$ is an obvious corollary of Theorem~\ref{t1}. Let
$\varepsilon> 0$ be as in the hypothesis of Theorem~\ref{t5} and let
$\tilde{w}=(w_0,w_1,0)$. A straightforward argument yields
\begin{align*}
   D \cap \spec\!_{\varepsilon,\tilde{w}}\, P_n& =\{z\in D:
   \sigma_n(z)\leq \varepsilon (w_0+w_1|z|)\}\\ &
   \subset \Big\{z\in D:
   \sigma_n(z)<\frac{\delta^2\gamma}{2(\delta^2+b_4)}\Big\}.
\end{align*}
Thus, if $n>N_3$, by virtue of Lemma~\ref{t8}$,
   C\cap\spec\!_{\varepsilon,\tilde{w}}\, P_n
    =\varnothing
$ guaranteeing the first conclusion of the claimed result.

On the other hand, Lemma~\ref{t7} ensures the existence of $N_4\geq
N_1$ such that
$\sigma_n(\lambda)<\frac{\delta^2\gamma}{2(\delta^2+b_4)}$ for all
$n>N_4$. Let $N:=\max\{N_3,N_4\}$ and assume that $n>N$. By
Lemma~\ref{t2}, there exists $z_n\in \spec P_n$ such that
$|z_n-\lambda|<\delta$. Furthermore, the connected component
$\Omega_n\subset \spec\!_{\varepsilon,\tilde{w}}\, P_n$ such that
$z_n\in \Omega_n$, satisfies $\Omega_n\subset D\setminus C$. Hence,
the second and third conclusions of the theorem follow directly from
Lemma~\ref{t4}.

\section{The essential spectrum}

The essential spectrum of $M$ is usually found by means of
analytical methods. Nonetheless, besides of being a natural question
\emph{per se}, the numerical evidence in all the examples discussed
in \cite{bou1}, \cite{lesh}, \cite{bou2} and \cite{bou3}, suggest
that approximation of this portion of the spectrum may also occur in
the ``second-order'' projection method described above. In this
final section, we discuss some results and open questions related to
this issue.

Ideally, we would like to know where does the whole set $\spec P_n$
accumulates in the limit $n\to \infty$. To this end, we may consider
the following two limiting set and then study the connection between
them as well as their relationship with $\spec M$. Given
$\varepsilon\geq 0$, let
\begin{gather*}
  \Lambda_\varepsilon:=\{\zeta \in \C: \mathrm{there\ exists\ }
  z_n\to \zeta\ \mathrm{such\ that}\
  \sigma_n(z_n)\leq \varepsilon \} \\
  \mathrm{and} \qquad
  \Sigma_\varepsilon:=\{ \zeta \in
  \C:\mathrm{limsup}_{n\to \infty}\,\sigma_n(\zeta )\leq \varepsilon \}.
\end{gather*}
Theorem~\ref{t1} ensures that $\spec\!_\mathrm{disc}\, M\subseteq
\R \cap \Lambda_0$. We now ask, which conditions yield
a similar property for $\spec\!_\mathrm{ess}\, M$. We do not intend
to answer this question here. Nevertheless,
elementary properties of these set might provide an insight towards
further investigations in this direction. Notice that $\Lambda_\varepsilon$
is the (uniform) limit set of $\spec\!_{\varepsilon,(1,0,0)} P_n$ as
$n\to \infty$.

For simplicity we assume from now on that $M$ is bounded. Thus
\begin{equation} \label{e9}
\Lambda_0\subseteq \Sigma_0 =\bigcap_{\delta>0} \Sigma_\delta =
\bigcap_{\delta>0} \Lambda_\delta.
\end{equation}
Indeed, directly from the definition, it follows that
\begin{equation} \label{e26}
\begin{aligned}
  \sigma_n(z)&\leq \min_{v\not = 0}\frac{\|P_n(w)v\|+
  |z-w|\|\Pi_n(2M-z-w)\Pi_nv\|}{\|v\|} \\
  &\leq \sigma_n(w)+|w-z| \sup\frac{\|(2M_n-z-w)v\|}{\|v\|} \\
  &\leq \sigma_n(w)+|w-z|c
\end{aligned}
\end{equation}
for all $|w-z|<\tilde{\varepsilon}<1$, where $c>0$ is chosen
independently from $\tilde{\varepsilon},n,z$ and $w$, because of $M$
is bounded. Now, assume on the one hand that $z\not
\in\Sigma_\varepsilon$. Then there exists $n(j)\in \N$ and $a>0$,
such that $\sigma_{n(j)}(z)>\varepsilon +a$ for all $j\in \N$.
According to \eqref{e26},
\[
   \sigma_{n(j)}(w)\geq \sigma_{n(j)}(z)-c|w-z|>\varepsilon+a-
   c\tilde{\varepsilon}
\]
for all $|w-z|<\tilde{\varepsilon}$. Hence, by choosing
$\tilde{\varepsilon}=a/(2c)$, it becomes evident that
$\sigma_{n(j)}(w)>\varepsilon$ whenever $|w-z|<\tilde{\varepsilon}$
for all $j\in \N$, so that $z\not \in \Lambda_\varepsilon$. Thus
\begin{equation*}
\Lambda_\varepsilon \subseteq \Sigma_\varepsilon \qquad\qquad
\mathrm{for\ all}\ \varepsilon\geq 0.
\end{equation*}
On the other hand, it is not difficult to prove that
\begin{equation*}
\Sigma_\varepsilon \subseteq \Lambda _{\varepsilon+\delta} \qquad
\mathrm{for\ all}\ \varepsilon\geq 0\ \mathrm{and}\ \delta>0
\end{equation*}
and that $\bigcap_{\delta>0}\Sigma_\delta \subseteq \Sigma_0.$ These
three inclusions ensure \eqref{e9}.

\medskip

The following proposition is crucial in the method suggested recently  by
Davies and Plum in \cite{sp}.

\begin{proposition}
If $M$ is bounded, then $\Sigma_0\,\cap\, \R=\spec M$.
\end{proposition}
\proof Indeed, if $\lambda\in \spec M$, for each $k>0$ there is
$\psi_k\in \CH$, $\|\psi_k\|=1$, such that
$\|(\lambda-M)^2\psi_k\|<1/k$. Then
\[
   \lim_{n\to\infty}
   \frac{\|P_n(\lambda)\Pi_n\psi_k\|}{\|\Pi_n\psi_k\|}< 1/k
\]
and so $\sigma_n(\lambda)\to 0$ as $n\to\infty$. Conversely, notice
that if $\lambda\in \R$ but $\lambda\not \in \spec M$, then
$
  \mathrm{Num}\, P_n(\lambda)\subset
  \mathrm{Num}\, (\lambda-M)^2 \subset [\mu,\infty).
$ Hence
\[
  \sigma_n(\lambda) \geq \dist[0,\mathrm{Num}\,P_n(\lambda)]
   \geq \mu,
\]
so that $\lambda \not \in \Sigma_0 \cap \R$.

\medskip

In other words, the inclusion $\Sigma_0\subseteq \Lambda_0$
complementary to the first one in \eqref{e9}, will automatically
imply approximation to the whole spectrum, in particular to
$\spec\!_\mathrm{ess}\, M$. For instance,
\cite[Proposition~3]{bou2}, if $M^2=Id$, then
$\Sigma_0=\Lambda_0\subseteq \{|\zeta|=1\}$. The validity of this
inclusion is closely related to the problem of whether there is an
upper bound independent of $n$ for the size of the blocks in the
Jordan canonical form of $A_n$, in the factorization
$P_n(z)=(z-A_n)(z-A_n^\ast)$. Indeed, let
$R_n(z)=(z-S_n)(z-S_n^\ast)$, where
\[S_n=\begin{pmatrix} 0 & & &
\\ 1&0 &  &  \\ & \ddots & \ddots & \\ & & 1 &0 \end{pmatrix}\in
\C^{n\times n}.\] Then $\spec R_n=\spec S_n \cup \spec
S_n^\ast=\{0\}$ for all $n\in \N$. Thus $\Lambda_0=\{0\}$. By
choosing $v=(z^{n-1},\ldots,z,1)^t$, $R_n(z)v=(z^{n+1},0,\ldots,0)$
so $\sigma_{R_n}(z)\leq |z|^{n+1}$. Hence $\Sigma_0= \{|z|\leq 1\}$.
Of course, although strict inclusion holds in this case, it not
clear whether $R_n(z)=P_n(z)$ for some $M=M^\ast$ and $\{\Pi_n\}$.


\subsection*{Acknowledgement} The author wishes to express his
gratitude to E.B.~Davies, P.~Lancaster and E.~Shargorodsky, for very
helpful discussions during the preparation of this manuscript.


\vspace{.5in}

\begin{minipage}{3in}
{\scshape Lyonell Boulton}\\
{\footnotesize Department of Mathematics and
Statistics, \\ University of Calgary, \\
Calgary, AB, Canada T2N 1N4 \\
email: \texttt{lboulton@math.ucalgary.ca}}
\end{minipage}

\end{document}